\documentclass[12pt]{article}

\usepackage{e-jc}
\usepackage{graphicx}

\usepackage{amssymb,latexsym,amscd, amsthm, amsmath, enumerate, subfigure}
\newtheorem{theorem}{Theorem}[section]

\newtheorem{proposition}[theorem]{Proposition}

\theoremstyle{definition}
\newtheorem{definition}[theorem]{Definition}
\newtheorem{example}[theorem]{Example}

\theoremstyle{remark}

\newtheorem{question}[theorem]{Question}

\date{\dateline{\today}}

\begin{document}

\title{Dissimilarity vectors of trees are contained in the tropical
Grassmannian}
\author{Benjamin Iriarte Giraldo
\\Department of Mathematics, San Francisco State University\\
San Francisco, CA, USA\\ {\tt biriarte@sfsu.edu}}
\maketitle

\begin{abstract}
In this short writing, we prove that the set of $m$-dissimilarity vectors
of phylogenetic $n$-trees is contained in the tropical Grassmannian $\mathcal{G}_{m,n}$, 
answering a question of Pachter and Speyer.
We do this by proving an 
equivalent conjecture proposed by Cools.

\end{abstract}
\section{Introduction.}

This article deals with the connection between phylogenetic trees
and tropical geometry. That these two subjects are
mathematically related can be traced back to 
Pachter and Speyer~\cite{PachterSpeyer}, Speyer and Sturmfels~\cite{SturmfelsSpeyer}, and
Ardila and Klivans~\cite{ardila}. The precise nature of this connection has been
the matter of some recent papers by Bocci and Cools \cite{BocciCools}
and Cools \cite{Cools}. In particular, a relation between $m$-dissimilarity
vectors of phylogenetic $n$-trees with the tropical Grassmannians $\mathcal{G}_{m,n}$
has been noted. 

\begin{theorem}[Pachter and Sturmfels~\cite{SturmfelsPachter}]\label{th:2dis}
The set of $2$-dissimilarity vectors
is equal to the tropical Grassmannian $\mathcal{G}_{2,n}$.
\end{theorem} 
This naturally raises
the following question.  
\begin{question}[Pachter and Speyer~\cite{PachterSpeyer}, Problem 3]\label{q:laqes}
Does the space of $m$-dissimilarity vectors lie in $\mathcal{G}_{m,n}$
for $m\geq 3$?
\end{question}
The result in this article is of relevance in this
direction and it is based on two papers
of Cools~\cite{Cools} and Bocci and Cools~\cite{BocciCools}, where the cases $m=3$, $m=4$ and
$m=5$ are handled. We answer Question~\ref{q:laqes} affirmatively for all
$m$:  
\begin{theorem}\label{col:final}
The set of $m$-dissimilarity vectors of phylogenetic $n$-trees is
contained in the tropical Grassmannian $\mathcal{G}_{m,n}$. 
\end{theorem}
\par
As we said, we prove Theorem \ref{col:final} by
proving an equivalent conjecture, Proposition~\ref{th:thet}
of this paper, or see Conjecture 4.4 of \cite{Cools}. 

\section{Definitions.}
\subsection{The Tropical Grassmannian.}

Let
$\mathbb{K}=\mathbb{C}\{\{t\}\}$ be the field of Puiseux
series. Recall that this
is the algebraically closed field of formal expressions
$$\omega=\sum_{k=p}^{\infty}c_kt^{k / q}$$
where $p\in\mathbb{Z}$, $c_p\neq 0$, $q\in\mathbb{Z}^+$
and $c_k\in\mathbb{C}$ for all $k\geq p$.
It is the algebraic closure of the field of Laurent series over
$\mathbb{C}$. The field comes equipped with a standard
valuation val$:\mathbb{K}\mapsto\mathbb{Q}\cup\{\infty\}$ by which
val$(\omega)=p/q$. As a convention, val$(0)=\infty$. 
\par
Now, let $x=(x_{ij})$ be an $m\times n$ matrix of indeterminates and let
$\mathbb{K}[x]$ denote the polynomial ring over $\mathbb{K}$ generated by these 
indeterminates.
Fix a second polynomial ring in $\binom{n}{m}$ indeterminates over
the same field: 
$$\mathbb{K}[p]=\mathbb{K}[p_{i_1,i_2,\dots,i_m}:1\leq i_1<i_2<\dots
<i_m\leq n]$$ 
Let $\phi_{m,n}:
\mathbb{K}[p]\mapsto \mathbb{K}[x]$ be the homomorphism of rings taking
$p_{i_1,\dots,i_m}$ to the maximal minor of $x$ obtained 
from columns $i_1,\dots,i_m$. 
\begin{definition}
The {\em Pl\"ucker ideal} or ideal of {\em Pl\"ucker relations} 
is the homogeneous prime ideal $I_{m,n}=$ker$(\phi_{m,n})$ which consists
of the algebraic relations or syzygies among the $m\times m$ minors of any
$m\times n$ matrix with entries in $\mathbb{K}$.  
\end{definition} 
For $m\geq 3$, the Pl\"ucker ideal has a Gr\"obner basis
consisting of quadrics; a comprehensive study of these ideals
can be found in Chapter 14 of the book by Miller and Sturmfels 
\cite{SturmfelsMiller} and in Sturmfels~\cite{Sturmfels}. It is a polynomial ideal in $\mathbb{K}[p]$ and we can define its {\em tropical variety}
in the usual way as we now recall. Let $a={\binom{n}{m}}$ and 
$\overline{\mathbb{R}}=\mathbb{R}\cup\{\infty\}$. Consider
$$f=\displaystyle\sum c_{\alpha}p^{\alpha_1}_{\sigma_1}p^{\alpha_2}_{\sigma_2}\dots
p^{\alpha_{a}}_{\sigma_{a}}\in\mathbb{K}[p],\mbox{ where
$\sigma_1,\dots,\sigma_a$ are the $a$ $m$-subsets 
of $\{1,\dots,n\}$}$$
\par
The {\em tropicalization} of $f$ is given by
$$\mbox{trop}(f)=\displaystyle\min\{
\mbox{val}(c_\alpha)+\alpha_1p_{\sigma_1}+\alpha_2p_{\sigma_2}+\dots+
\alpha_{a}p_{\sigma_a}\}.$$ 
The {\em tropical
hypersurface} $\mathcal{T}(f)$ of $f$ is the   
set of points in $\overline{\mathbb{R}}^a$ where
trop$(f)$ attains its minimum twice or, equivalently, 
where trop$(f)$
is not differentiable. 
\par
We are now ready to define tropical Grassmannians.
\begin{definition}
The tropical variety 
$\mathcal{T}(I_{m,n})=\displaystyle\bigcap_{f\in I_{m,n}}\mathcal{T}(f)$
of the Pl\"ucker ideal $I_{m,n}$ is denoted by $\mathcal{G}_{m,n}$
and is called a {\em tropical Grassmannian}. 
\end{definition} 
\par
We have the following fundamental characterization of $\mathcal{G}_{m,n}$ 
which is a direct application of a more general 
result~\cite[Theorem 2.1]{SturmfelsSpeyer}.
\begin{theorem}
The following subsets of $\overline{\mathbb{R}}^{a}$
coincide:
\begin{itemize}
    \item The tropical Grassmannian $\mathcal{G}_{m,n}$.
    \item The closure of the set $\{(val(c_1),val(c_2)
,\dots,val(c_a)):(c_1,c_2,\dots,c_a)\in 
V(I_{m,n})\subseteq\mathbb{K}^{a}\}$
\end{itemize}
\end{theorem}
\subsection{\textsf{Phylogenetic Trees.}}
We also treat phylogenetic trees in this paper. 

\begin{definition}
A {\em phylogenetic $n$-tree } 
is a tree which has a labeling of its $n$ leaves
with the set $\{1,\dots,n\}$ and such that each
edge $e$ has a positive
real number $w(e)$ associated to it, which we call the {\em weight}
of $e$.  
\end{definition}
There is also a crucial related family of trees which we now define:
\begin{definition} An {\em ultrametric $n$-tree} is 
a binary rooted tree which has a labeling of its $n$ leaves
with $\{1,\dots,n\}$ and such that
\begin{list}{$\bullet$}{\setlength\itemsep{-0.05in}}
\item each
edge $e$ has a nonnegative
real number $w(e)$ associated to it, called the {\em weight}
of $e$
\item it is $d$-equidistant, for some $d>0$, {\it i.e.} the sum of
the edges in the path from the root to every leaf is
precisely $d$ 
\item the sum of the weights of all edges in the path connecting every two different leaves is positive.
\end{list}
\end{definition}
\par
Particularly, note 
that an ultrametric tree is binary and may have edges of weight $0$.
Now, let $T$ be a phylogenetic $n$-tree. 
Define the vector $D(m,T)$ whose entries
are the numbers $d_\sigma$, where
$\sigma$ is a subset of $\{1,2,\dots,n\}$ of size $m$ and
$d_{\sigma}$ is the {\em total weight} of the smallest subtree
of $T$ which contains the leaves in $\sigma$. By the total weight
of a tree, we mean the sum
of the weights of all the edges in that tree.  

\begin{definition}
The vector $D(m,T)$ is called the {\em $m$-dissimilarity
vector of $T$}. The set of all $m$-dissimilarity vectors
of phylogenetic trees with $n$ leaves will be called
the {\em space of $m$-dissimilarity vectors of $n$-trees}.
\end{definition}  
\begin{definition}
A metric space $S$ with distance function 
$d:S\times S\mapsto \mathbb{R}_{\geq 0}$ 
is called an {\em ultrametric space} if the
following inequality holds for all $x,y,z\in S$:
$$d(x,z)\leq \mbox{max}\{d(x,y),d(y,z)\}$$
\end{definition}
It is a well known fact that 
finite ultrametric spaces are 
realized by ultrametric trees, see for example~\cite[Lemma 11.1]{BB}.  
\subsection{Column Reductions.}

Let $n\geq 4$. Suppose we are given integers  
$1\leq a,b\leq n$ with $a\neq b$ and let
$c_{a,b}$ be the operator acting on Puiseux matrices for which, 
for any $n\times n$ matrix $M$, $c_{a,b}(M)$ is the matrix 
obtained from $M$ by subtracting column $b$ to column $a$.
We know $c_{a,b}$ preserves the determinant, {\it i.e.} 
$\mbox{det}\left(c_{a,b}(M)\right)=\mbox{det}(M)$.  
For $l\geq 1$, let 
$\left(c_{a_l,b_l}\circ\dots\circ c_{a_2,b_2}\circ c_{a_1,b_1}\right)(M)$ be the matrix
obtained from $M$ by first subtracting 
column $b_1$ to column $a_1$, then subtracting
column $b_2$ to column $a_2$, and so on up to
subtracting column $b_l$ to column $a_l$. Call
this matrix a {\em column
reduction of $M$} if 
the following conditions are met:
\begin{itemize}
\item $1\leq a_1,\dots,a_l,b_1,\dots,b_l\leq n$
\item the numbers $a_1$, $a_2$, $\dots$, $a_l$ are pairwise different
\item whenever $1\leq k\leq l$, the number
$b_k$ is different from $a_{1}$, $\dots$, $a_k$.
\end{itemize}
For simplicity, we will
accept $M$ as a column reduction
of itself.  

\section{Main Result.}\label{sec:mr}

We are now ready to prove 
Theorem~\ref{col:final}. Cools~\cite{Cools} reduced it
to the following statement which we now prove.  

\begin{proposition}[Cools~\cite{Cools}, Conjecture 4.4 ]\label{th:thet}
Assume $n\geq 4$. Let $T$ be a $d$-equidistant ultrametric $n$-tree
with root $r$ and such that all its edges have rational
weight. 

For each edge $e$ of $T$, denote by $h(e)$ the well-defined 
sum of the weights of all the 
edges in the path from the top node of $e$ to any leaf below $e$ and
let $a_1(e),\ldots,a_{n-2}(e)$ be generic complex numbers. 
 
Let $x_i^{(j)}\in\mathbb{K}$ 
(with $i\in\{1,\ldots,n\}$ and $j\in\{1,\ldots,n-2\}$) be the sum of the monomials $a_j(e)t^{-h(e)}$, where $e$ runs over all edges between $r$ and $i$. 
Then, the valuation of the determinant 
of 
$$M=\begin{pmatrix} 1&1&\ldots&1\\ x_1^{(1)}&x_2^{(1)}&\ldots&x_n^{(1)}\\ (x_1^{(1)})^2&(x_2^{(1)})^2&\ldots&(x_n^{(1)})^2\\ x_1^{(2)}&x_2^{(2)}&\ldots&x_n^{(2)} \\ \vdots&\vdots&\vdots&\vdots \\ x_1^{(n-2)}&x_2^{(n-2)}&\ldots&x_n^{(n-2)}\end{pmatrix}$$
is equal to
$-D$, where $D$ is the total
weight of $T$.
\end{proposition}
In the course of the proof,
we assume $T$ is binary, which follows from
the construction of Bocci and Cools~\cite{BocciCools}.
Notice they start with a phylogenetic tree and then
define an associated ultrametric from its $2$-dissimilarity vector, 
therefore inducing an ultrametric tree. Here, $T$ 
corresponds to certain subtrees of this induced ultrametric tree. 
\begin{proof}

As $T$ is binary, we know $T$ has $n$ leaves, 
$n-2$ internal nodes of degree
three, one node (the root) of degree two and
$2(n-1)$ edges. \\ 
\par
Let $\leq_T$ be the tree order of $T$ with respect to
$r$, {\it i.e.} the order on the set of nodes of $T$ by which
$v\leq_T w$ iff $v$ lies in the path from $r$ to $w$ in $T$. Let $v_1,v_2,\dots,v_{n-1}$ be the $n-1$ internal nodes of $T$
numbered in such way that if $v_i\leq_T v_j$, then $j\leq i$.
We must have $v_{n-1}=r$.  \\
\par
Consider an injective function $\alpha:v_i\mapsto a_i$
from the set of internal nodes to the leaves of $T$
so that $v_i\leq_T a_i$ for all $i$ with $1\leq i\leq n-1$. 
Now, for each of these values of $i$, let
$b_i$ be the unique leaf such that $b_i\neq a_j$ for all $j$ with 
$1\leq j\leq i$, and such that $v_i\leq_T b_i$. \\
\par
To show the existence of $\alpha$, we construct it 
succesively starting with
$\alpha(v_1)$, then $\alpha(v_2)$ and then continuing up until
we define $\alpha(v_{n-1})$. Suppose we have already defined
$\alpha(v_1)$, $\dots$, $\alpha(v_{i-1})$ for
some $i<n-1$. Consider the maximal subtree $T_i$ of $T$ whose root is
$v_{i}$, {\it i.e.} $T_i$ is the subtree below $v_{i}$. If this tree has
$m$ leaves then it has $m-1$ internal nodes, including $v_{i}$ itself.
So far we haven't defined $\alpha$ for nodes
between $r$ and $v_i$ but we have defined
it for all internal nodes of $T_i$ different from
$v_i$. Therefore, there are exactly
$m-2$ leaves of the tree $T_i$ which have been assigned
to some of $v_1$, $\dots$, $v_{i-1}$ under $\alpha$, so there are $2$ leaves
which we can assign to $v_i$: $\alpha(v_i)$ can be either one
of them. Incidentally this also gives us the existence and uniqueness
of the respective $b_i$. \\
\par
Now, we want to establish the equality 
$\sum_{i=1}^{n-1}h(v_i)=D-d$. This equality
is clearly true when $T$ has $2$ or $3$ leaves, so that $n=2$
or $n=3$. Let
now $n>4$ and 
suppose we have proved the result for all trees with $i$ leaves with 
$i<n$. Recall $n$ is being taken as the number of leaves in $T$, which
is rooted $d$-equidistant with root $r=v_{n-1}$. 
We know the equality holds for each of the subtrees
$T_1$, $\dots$, $T_{n-2}$ below
$v_{1}$, $\dots$, $v_{n-2}$, respectively. Let $T_{n-2}$ be
$d_{n-2}$-equidistant and let $T_{n-3}$ be $d_{n-3}$-equidistant.  
There are two cases to distinguish.
If $v_{n-2}<_T v_{i}$ for all $i<n-2$ then
$\sum_{i=1}^{n-2}h(v_i)=\left(D-d-(d-d_{n-2})\right)-
d_{n-2}=D-2d$ by induction, so $\sum_{i=1}^{n-1}h(v_i)=D-d$. 
Otherwise suppose $v_{n-2}$ and $v_{n-3}$ are incomparable
in $<_T$. 
Then $T_{n-2}$ and $T_{n-3}$ are disjoint graphs and we have
$$\displaystyle\sum_{v_i\in V(T_{n-2})}h(v_i)=
\left(D-\left(\displaystyle\sum_{v_j\in V(T_{n-3})}h(v_j)+d_{n-3}\right)-
(d-d_{n-3})-(d-d_{n-2})\right)-d_{n-2}$$   
by induction. Reordering we get
$$\displaystyle\sum_{v_i\in V(T_{n-2})}h(v_i)+
\displaystyle\sum_{v_j\in V(T_{n-3})}h(v_j)=D-2d$$ 
so if we add $h(v_{n-1})=d$ to both sides we get our result. \\
\par
Now consider the column reduction $M^* =\left(c_{a_{n-1},b_{n-1}}\circ 
\dots\circ c_{a_2,b_2}\circ 
c_{a_1,b_1}\right)(M)$ of $M$. We claim that the valuation of all nonzero
monomials
$\prod_{i=1}^{n}M^{*}_{i,\sigma
(i)}$ with $\sigma\in S_n$ in the sum
$$\mbox{det}(M^*)=
\sum_{\sigma\in S_n}\left(\mbox{sgn}(\sigma)\prod_{i=1}^{n}M^{*}_{i,\sigma
(i)}\right),$$
is precisely $-\left(\sum_{i=1}^{n-1}h(v_i)+d\right)=-D$.
To see this notice for all $i$, $1\leq i\leq n-1$, we have 
\begin{itemize}
\item $M^{*}_{1a_i}=0$
\item the valuation of
$M^{*}_{3a_i}$ is $-d-h(v_i)$
\item the valuation of $M^{*}_{ja_i}$ is $-h(v_i)$ if
$j\neq 1$ and $j\neq 3$
\item the only nonzero term in the first row of $M^{*}$ is the $1$ in column
$b_{n-1}$
\end{itemize}

Because of our generic choice of coefficients, we can find some
monomial term in the sum $\mbox{det}(M^*)$ with valuation $-D$ 
which doesn't get cancelled, so we are done. \\  
\end{proof}
\begin{figure}[ht]
\centering
\includegraphics{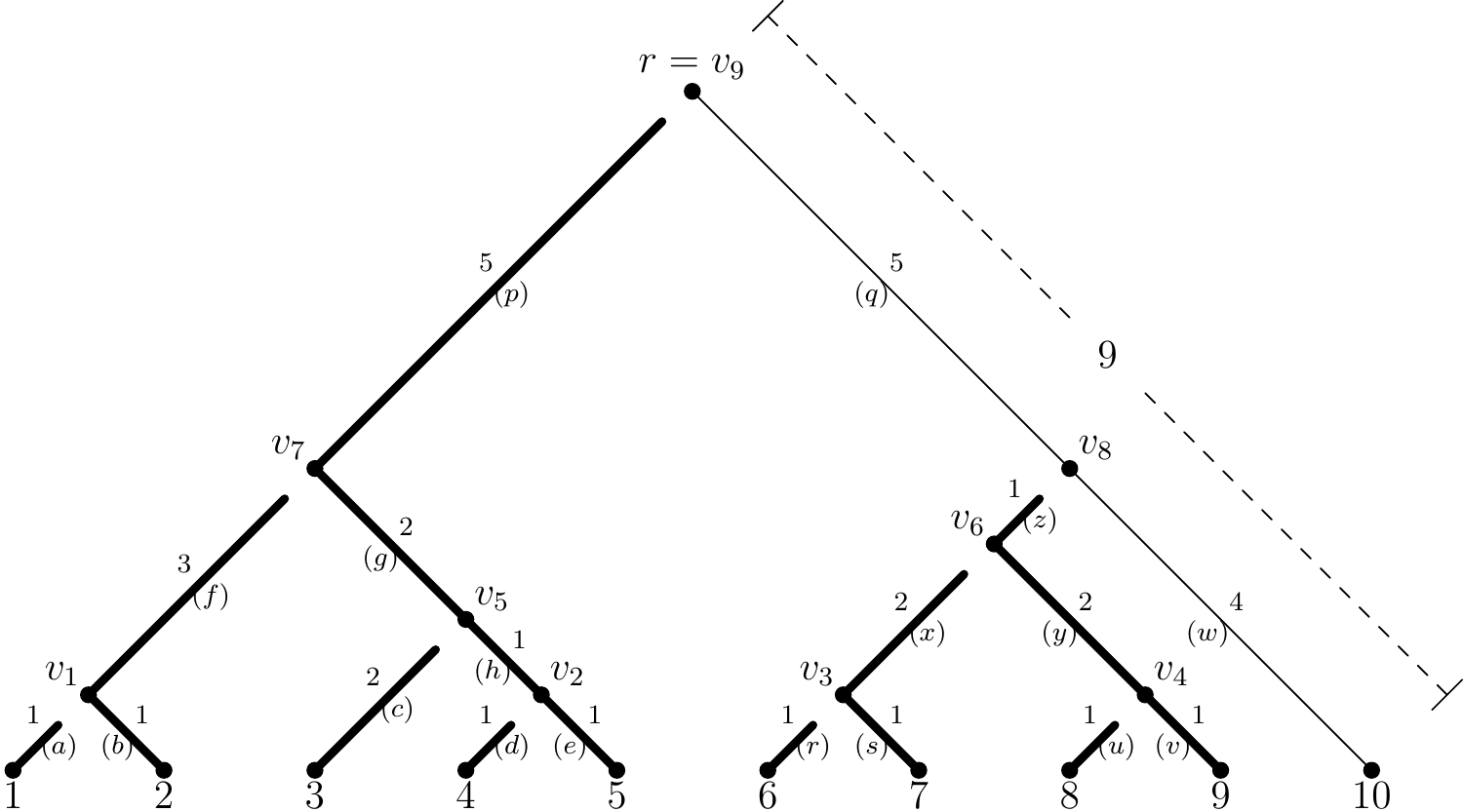}
\caption{A rooted $10$-tree. The injective function
$\alpha:=\{(v_1,1),(v_2,4),(v_3,6),(v_4,8),(v_5,3),(v_6,7),
(v_7,2),(v_8,9),(v_9,5)\}$ is depicted, as well as the
equality $\sum_{i=1}^9 h(v_i)=35-9$.}
\label{fig:fourt}
\end{figure}
\begin{example}
Consider the $9$-equidistant $10$-tree of Figure~\ref{fig:fourt} with
total weight $35$. The second row
of the matrix $M$ associated to this tree is the following vector
with generic complex coefficients:
\begin{align*}
\nonumber [ & at^{-1}+ft^{-4}+pt^{-9} & , & bt^{-1}+ft^{-4}+pt^{-9}  &
,    &
ct^{-2}+gt^{-4}+pt^{-9} &, \\
\nonumber & dt^{-1}+ht^{-2}+gt^{-4}+pt^{-9}  & ,  & 
et^{-1}+ht^{-2}+gt^{-4}+pt^{-9} & ,  &
rt^{-1}+xt^{-3}+zt^{-4}+qt^{-9} &, \\
\nonumber & st^{-1}+xt^{-3}+zt^{-4}+qt^{-9} & ,
& ut^{-1}+yt^{-3}+zt^{-4}+qt^{-9} & 
, &  
vt^{-1}+yt^{-3}+zt^{-4}+qt^{-9} &, \\
\nonumber & wt^{-4}+qt^{-9}] 
\end{align*}
Using the operator
$\left(c_{5,10}\circ c_{9,10}\circ c_{2,5}\circ c_{7,9}\circ
c_{3,5}\circ c_{8,9}\circ c_{6,7}\circ c_{4,5}\circ c_{1,2}\right)$
suggested by the figure we obtain the column reduction 
$M^*$ whose second row is the vector: 
\begin{align}
\nonumber [ & (a-b)t^{-1}&, & \ \ \ \ \ (b-e)t^{-1}-ht^{-2}+(f-g)t^{-4} &,\\
\nonumber & -et+(c-h)t^{-2} &, & \ \ \ \ \ (d-e)t^{-1} &, \\
\nonumber & et^{-1}+ht^{-2}+(g-w)t^{-4}+(p-q)t^{-9} &
,& \ \ \ \ \ (r-s)t^{-1}&, \\
\nonumber & (s-v)t^{-1}+(x-y)t^{-3}&,& \ \ \ \ \ (u-v)t^{-1} &, \\
\nonumber & vt^{-1}+yt^{-3}+(z-w)t^{-4} &, & \ \ \ \ \ wt^{-4}+qt^{-9}] 
\end{align}
It has valuation vector: 
\begin{align}
\nonumber (&-1,-4,-2,-1,-9,-1,-3,-1,-4,-9) = & \\
\nonumber (&-h(v_1),-h(v_7),-h(v_5),-h(v_2),-h(v_9),-h(v_3),-h(v_6),
-h(v_4),-h(v_8)) & ,
\end{align}
where $v_1,v_7,v_5,v_2,v_9,v_3,v_6,v_4,v_8$ are the preimages of
$1,2,3,4,5,6,7,8,9$ under $\alpha$, respectively in that order. 
Also notice that $\sum_{i=1}^{9}h(v_i)=35-9$.
\end{example}
We have shown that the $m$-dissimilarity vector of a phylogenetic tree 
$T$ with $n$ leaves gives a point in the tropical Grassmannian
$\mathcal{G}_{m,n}$, and therefore gives rise to 
a tropical linear space. 
The combinatorial structure of those tropical linear
spaces is the subject of an upcoming paper~\cite{Iriarte}.

\section{Acknowledgements.}
This work began to develop itself at Federico Ardila's course
on Combinatorial Commutative Algebra, jointly
offered at San Francisco State University and the
Universidad de los Andes in the spring
of 2009. Special thanks
to Federico for many useful commentaries and suggestions, including
a beautiful simplification of my original proof of 
Lemma~\ref{th:thet} and for
bringing to my knowledge the paper of Cools
~\cite{Cools} and Question~\ref{q:laqes}. 
Thanks to the SFSU-Colombia Combinatorics Initiative for supporting
this research project.  
\bibliography{finalproject}{}
\bibliographystyle{amsplain}
\providecommand{\bysame}{\leavevmode\hbox to3em{\hrulefill}\thinspace}
\providecommand{\MR}{\relax\ifhmode\unskip\space\fi MR }
\providecommand{\MRhref}[2]{%
  \href{http://www.ams.org/mathscinet-getitem?mr=#1}{#2}
}
\providecommand{\href}[2]{#2}

\end{document}